\documentclass[14pt]{article}
\usepackage{amssymb}
\usepackage{amsmath}
\usepackage{latexsym}
\usepackage{epsfig}
\newcommand {\eq}[1]{\begin{equation}\label{#1}}
\newcommand {\en} {\end{equation}}
\newcommand {\proof} {\noindent{\it Proof}. \ignorespaces}
\newcommand {\eproof}{\space{\ \vbox{\hrule\hbox{\vrule height1.3ex\hskip0.8ex\vrule}\hrule}}\vskip 0.3cm \par}
\newcommand {\R}        {{\mathbb R}}

\newcommand {\Rn}       {\R^n}

\newcommand {\Rp}       {\R_+}
\newcommand {\Rpn}       {\R^n_+}
\newcommand {\Rpk}       {\R^k_+}

\newcommand {\Rpnk}      {\R^{nk}_+}
\newcommand {\mat}[1]{\left[\begin{array}{#1}}
\newcommand {\rix}          {\end{array}\right]}

\newtheorem{theorem}          {Theorem}
\newtheorem{lemma}         [theorem]{Lemma}
\newtheorem{definition}    [theorem]{Definition}
\newtheorem{corollary}     [theorem]{Corollary}
\newtheorem{proposition}     [theorem]{Proposition}
\newtheorem{key}    [theorem]{Key Proposition}

\newtheorem{example}        [theorem]   {Example}
\newtheorem{remark}         [theorem]   {Remark}

\newtheorem{algorithm}      [theorem]  {Algorithm}

\def\C{{\rm C\kern-.48em\vrule width.06em height.6em depth-.02em
        \kern.48em}}
\def\gl{\lambda}

\def\ge{\varepsilon}

\def\msn{\medskip\noindent}

\def\nn{\{1,\ldots,n\}}
\newcommand {\supp}     {\mathop{\rm supp}\nolimits}
\newcommand {\spn}    {\mathop{\rm span}\nolimits}
\begin{document}
\title{Generators, Extremals and Bases of Max Cones\thanks{%
Supported by the EPSRC grant T06414 and the RFBR grant 05-01-00824}}
\author{Peter Butkovi\v{c}\thanks{Corresponding author.}\;\thanks{%
School of Mathematics, University of Birmingham, Edgbaston,
Birmingham B15 2TT, UK.} \and {Hans Schneider\thanks%
{Department of
Mathematics, University of Wisconsin--Madison, Madison, WI 53706,
USA}}\and {Serge\u{\i} Sergeev\thanks%
{Department of Physics, Sub-Department of Quantum Statistics and
Field Theory, Moscow State University, Moscow,
119992 Vorobyevy Gory, Russia}}}%
 \maketitle
%{\small\version}
%\date{}
\begin{abstract}
\msn
Max cones are max-algebraic analogs of convex cones. In the present
paper we develop a theory of generating sets and extremals of max
cones in $\Rpn$. This theory is based on
the observation that extremals are minimal elements of max cones
under suitable scalings of vectors. We give new proofs of existing
results suitably generalizing, restating and refining them. Of
these, it is important that any set of generators may be partitioned
into the set of extremals and the set of redundant elements. We
include results on properties of open and closed cones, on
properties of totally dependent sets and on computational bounds for
the problem of finding the (essentially unique) basis of a finitely
generated cone.

\msn AMS classification: 15A48, 15A03.

\msn Keywords: Max algebra, cone, basis, scaling, extremal, algorithm.
\end{abstract}
\section{Introduction}%
By {\em max algebra} we understand the analog of linear algebra
obtained by considering $\R_+$ (the nonnegative reals)
with max times operations:
$$a\oplus b := \max(a,b),$$
$$a\otimes b := ab$$
extended to matrices and vectors. That is, if
$A=(a_{ij}),~B=(b_{ij})$ and $C=(c_{ij})$ are matrices of
compatible sizes with entries from $\R_+$ and $\alpha\in\R_+$, we
write $C=A\oplus B$ if $c_{ij}=a_{ij}\oplus b_{ij}$ for all $i,j$,
$C=A\otimes B$ if $c_{ij}=\bigoplus_{k} a_{ik}\otimes
b_{kj}=\max_{k}(a_{ik}b_{kj})$ for all $i,j$ and $C=\alpha \otimes
A$ if $c_{ij}=\alpha \otimes a_{ij} $ for all $i,j$. There are
several essentially equivalent \footnote{that is, algebraically
isomorphic} definitions of max algebra. An early paper presenting
the above version is \cite{Vor}, another early paper presenting an
equivalent version is \cite{CG0}.
%The algebraic structures appearing in tropical geometry, for instance in
%\cite{DS} and \cite{Jos}, are also to some extent
%another form of max algebra.
For more information on max algebra, its generalizations and applications
the reader is referred e.g. to \cite{BCOQ}, \cite{But}, \cite{CG1,CG2,CG3},
\cite{GM} and \cite{UZ}.
%Other forms
%of max algebra may be found e.g. in \cite{CG0,CG1,CG2,CG3} and in \cite{But}.
See also \cite{LM} for recent developments in the area and for
further references.

\msn
%As observed in the abstract, our aim is to develop a theory of
%general max cones in max algebra and here
We give a summary of the
contents of this paper. In Section 2 we begin by defining (max)
cones, extremals, generating sets, independent sets and totally
dependent sets. The key observation is Proposition \ref{spanit}
which extends \cite[Proposition 2.9]{Jos}. We deduce in Theorem
\ref{morextr} that extremals are minimal elements of max cones under
suitable scalings of vectors. This leads us to a reformulation and
new proof of the basic result Theorem \ref{mainres} which is also
easily derived from \cite[Proposition 2.5.3]{Wag1}: Every generating
set of a max cone can be partitioned into the set of the extremals
of a cone and a set of redundant elements. It follows that if a cone
has a basis then (under a scaling condition) it consists of the
extremals of the cone and hence the cone has a basis unique up to
scalar multiples, see Theorem \ref{nextres} and its corollaries for
more detail. In Corollary \ref{totdep} we then turn to the case when
the set of extremals of the cone is empty, in which case every
generating set is totally dependent. Discussion of totally dependent
sets specifically may be new. Towards the end of this section we
consider topological notions. In Corollary \ref{open} we show that
(under a natural restriction) every open cone has totally
dependent generators and in Proposition \ref{c:cc} we prove an
analogue of Minkowski's theorem for closed cones. This result
extends a result due to \cite{Hel} and it also appears as
\cite[Theorem 3.1]{GK1} where a different proof is given.
%We give proofs of known results found in \cite{Wag1} and \cite{Jos}
%suitably restating and refining them. Further results include
%Minkowski's theorem, which extends the results of \cite{Hel}, and
%some properties of totally dependent sets.

\msn In Section 3 we give two simple versions of an algorithm, based
on \cite[Theorem 16.2]{CG2} for finding the (essentially unique)
basis of a finitely generated max cone and a MATLAB program which
implements one version. We also relate our problem to the classical
problem of finding maxima of a set of vectors which is described in
\cite{KLP} and in \cite[Section 3]{PS}, and give the bounds for
computational complexity.

\msn  We now relate the concepts and techniques of our proofs to
those in other publications. Most of our concepts appear in \cite{Hel},
\cite{Wag0} and \cite{Wag1}, sometimes under different
names. For instance, extremals are called irreducible
elements in \cite{Wag0} and \cite{Wag1}, and minimal elements
are called efficient points in \cite{Hel}.
Our key Proposition \ref{spanit} may also be derived
from (possibly slightly extended) results found in some of our
references. Examples are results in \cite{Vor} Section 2 in terms
of set coverings, see also \cite{CG0}, \cite[Theorem 15.6]{CG2}
and \cite{But} Section 2, or the fundamental results of \cite{CG2}
and \cite[Section 3]{CG3} concerning max linear systems. The
latter are also found in \cite{CG1} and \cite{DS} in terms of a
projection operator.
Further, such results on max linear
systems as Proposition \ref{spanit} and Proposition \ref{extension}
can be extended to the case of
functional Galois connections, as it is shown in \cite{AGK}.
The generalizations considered
in \cite{AGK} are useful in many areas including abstract convex
analysis, the theory of Hamilton-Jacobi equation and the Monge-Kantorovitch
mass transportation problem.

\msn Our topic is also related to (and partially stimulated by) the
emerging field of tropical geometry which develops basic concepts of
max algebra in a different form and with different terminology and
applies these to finitely generated structures, see
\cite{DS},\cite{Jos} and \cite{BYu}. In particular, Proposition
\ref{spanit} can also be seen as a minor extension of
\cite[Proposition 2.9]{Jos} which is important in the theory of
tropical halfspaces. The emphasis of these papers is on geometry,
while in this paper it is on algebraic and order theoretic results.
\msn Max cones are also studied in \cite{GK0} and \cite{GK1}. The
main effort of these papers is to develop the theory of max-plus
convex sets and their recession cones. This theory is not present in
our paper. In turn, we deal with more general cones and we emphasize
the link to set maxima and give a more detailed description of bases
and generating sets.

\msn In max algebra as in linear algebra a basis is normally
defined as an {\em indexed set}, that is a {\em sequence} if the
basis is finite or countable, see \cite{Wag1} for a definition in
max algebra or \cite[p.10]{Bou} in linear algebra. Since we wish
to show the inclusion of the set of extremals (which do not have a
natural order) in every generating set or basis for a cone we
define the latter in term of {\em sets} in Section 2. We thereby
exclude the possibility of a repetition of elements in generating
sets. But we change our point of view in Section 3 on algorithms for
finitely generated cones since we wish to consider the generators
as columns of a matrix.

\msn Max cones have much in common with convex cones, see \cite{Roc}
for a general reference.  This has been exploited (and
generalized) in many papers including those just quoted and e.g.
\cite{CGQS} and \cite{KZ}. To this end, the basic concepts of this
paper and such results as Theorem \ref{mainres}, Proposition
\ref{c:cc} and Proposition \ref{p:clcl} have their direct analogs
in terms of positive linear combinations and
in convex analysis. We do not provide details, as convex geometry
is also beyond our scope here.
%\end{document} TEMP 55555
\section{Generating sets, bases and extremals}
\label{gsb&e}
\msn We begin with two standard definitions of max algebra.
\begin{definition}{\rm A subset $K$ of $\Rpn$ is a {\em max cone} in
$\Rpn$ if it
is closed under $\oplus$ and $\otimes$ by nonnegative reals.}
\end{definition}
\begin{definition} {\rm Let $S \subseteq \Rn_+$. Then $u$ is a
{\em max combination} of $S$ if
\begin{equation}\label{maxcomb}
u = \bigoplus_{x \in S} \gl_x x, \quad \gl_x \in \Rp,
\end{equation}
where only finite number of
%Sergei%
%nonnegative scalars
%%%%%%%
$\gl_x \neq 0$. The set of all max
combinations will be denoted by span($S$). We put span$(\emptyset)=\{0\}$.}
\end{definition}
Evidently, span($S$) is a cone. If span($S) = K $, we call $S$ a
{\em set of  generators} for $K$.
\begin{definition}\label{extr}  An element $u \in K$ is
an {\em extremal} in $K$ if
\begin{equation}
u = v \oplus w,\ v,w  \in K \Longrightarrow u = v\ {\rm or}\ u= w.
\end{equation}
\end{definition}
If $u$ is an extremal in $K$ and $\gl > 0$ then $\gl u$ is also an
extremal in $K$.
%%%%%%%%
%Extremals are defined in \cite{Wag1} under the name "irreducible
%elements".
%Sergei%
\begin{definition}{\rm
An element $x \in \Rn_+$ is {\em scaled} if $||x||=1$.}
\end{definition}
For most of this section, $||x||$ may be any norm in $\Rn$ (they are all
equivalent). However, in the end we specialize to the max norm,
$\|x\|=\max x_i$, in order to exploit the property that it is max
linear on $\Rpn$.
If $S \subseteq \Rpn$ we may call $S$ scaled to indicate that it
consists of scaled elements.
\begin{definition} {\rm Let $S$ be a set of vectors in $\Rpn$.
\begin{enumerate}
\item
The set $S$ is {\em dependent} if, for some $x \in S$, $x$  is a max combination
of $S\setminus\{x\}$. Otherwise,
$S$ is {\em independent}.
\item
The set  $S$ is {\em totally dependent} if every $x \in S$  is a
max combination of $S\setminus\{x\}$.
\end{enumerate}}
\end{definition}
Thus the empty set of vectors is both independent and totally
dependent. Since span$(\emptyset) = \{0\}$, the set $\{0\}$ is
totally dependent.
\begin{definition} \rm Let $K$ be a cone in $\Rpn$.
A set $S$ of vectors in $\Rpn$ is
a {\em  basis} for $K$ if it is an independent set of generators
for $K$.
\end{definition}
The set of all unit vectors $\{e^p,\ p=1,\ldots,n\}$
defined, as usual, by
\begin{equation}
\label{unit}
e_j^p=
\begin{cases}
1,& j=p\\
0,& j\ne p
\end{cases},
\end{equation}
is a basis of $\Rpn$, which is called {\em standard}.
\begin{lemma} \label{extr1} Let $S$ be a set of scaled generators
for the cone $K$ in $\Rpn$ and let $u$ be a scaled extremal in
$K$. Then $u \in S$.
\end{lemma}
\proof   Suppose $u$ is given by the max combination
(\ref{maxcomb}). Since the number of nonzero $\gl_x$ is finite, we
may use Definition \ref{extr} and induction to show that $u=\gl_x x$ for
some $x$. But $u$ and $x$ are both scaled, hence $u=x$ and $u \in S$.
\eproof
\begin{lemma}\label{extr2} The set of scaled extremals of a cone is
independent.  \end{lemma}
\proof If the set $E$ of scaled extremals is nonempty let $u$ be a
scaled extremal in $K$ and apply Lemma \ref{extr1} to the cone
$K_1 :=$ span$(E\backslash \{u\})$. This shows $u \not \in K_1$
and the result is proved. \eproof
Below we use subscripts for elements of vectors in $\Rpn$ and
superscripts to label vectors.
\begin{definition} {\rm
Let $v \in \Rpn$. Then the {\em support} of $v$ is defined by
$$ {\rm supp}(v) = \{j \in \nn: v_j > 0\}.$$
The cardinality of
$\supp(u)$ will be written as $|\supp(u)|$.}
\end{definition}
In order to relate the natural partial order on $\Rpn$ to results
on extremals of cones we introduce a scaling of vectors in $\Rpn$
for each $j \in \nn$ such that for each scaled vector $v_j = 1$.
\begin{definition}{\rm
\begin{enumerate}
\item
Let $u \in \Rpn$ and suppose $j \in \supp(u)$. Then we define $
u(j) = u/u_j.$
\item
Let $S \subseteq \Rpn$. We define $S(j) = \{u(j) : u \in S\ {\rm
and}\ j \in \supp(u)\}$ for all $j=1,\ldots,n$.
\item
Let $S \subseteq \Rpn$. An element $u \in S$ is called {\em
minimal} in $S$, if $v \leq u$ and $v \in S$ implies that $v = u$.
\item
Let $K$ be a cone in $\Rpn$, let $u \in K$,
%Sergei%
and let
$j \in \supp(u)$.
%%%%%%%%
We define
$$D_j(u) = \{v \in K(j): v \leq u(j)\}.$$
\end{enumerate}}
\end{definition}
Our key observation is the following Proposition. It can be viewed
as a minor but needed extension of \cite[Proposition 2.9]{Jos},
see also the remarks concerning it in our Introduction.
\begin{key} \label{spanit}
Let $S \subseteq \Rpn$. Then the following are equivalent:
\begin{enumerate}
\item
$u \in {\rm span}(S)$.
\item
For each $j \in \supp(u)$ there is an $x^j \in S$ such that $j \in
\supp(x^j)$  and $x^j(j) \in D_j(u)$.
\end{enumerate}
\end{key}
\begin{proof}
2. $\Longrightarrow$ 1: If 2. holds, then  $u =\bigoplus_{j
\in \supp(u)} \gl_j x^j$ where $\gl_j = u_j/x^j_j$.\\
1. $\Longrightarrow$ 2.
Conversely if 1. holds, then it follows immediately from (\ref{maxcomb}) that for
each $j \in \supp(u)$ there is an $x^j \in S$ with $\gl_j x^j \leq
u$ and $(\gl_j x^j)_j = u_j$. Clearly, $\gl_j = u_j/x^j_j$ which
yields 2.
\end{proof} \eproof
The following immediate corollary to Proposition \ref{spanit} is
essentially found as \cite[Theorem II.1]{Hel} and as
\cite[Proposition 5]{DS}. It is analogous to Carath\'{e}odory's Theorem.
\begin{corollary}\label{cara}
Let $S \subseteq \Rpn$. Then $u \in \spn(S)$ if and only if there
are k vectors $x^1,\ldots, x^k \in S,$ where  $k \leq |supp(u)|$,
such that $u \in \spn\{x^1,...,x^k$\}. \eproof
\end{corollary}
\begin{corollary}\label{remove} Let $K$ be a cone in $\Rpn$ and let $T$ be a set
of generators for $K$. Let $U \subseteq T$ and let  $S =
T\backslash U$. Then $S$ generates $K$ if and only if each $u \in
T$ satisfies condition 2. of Proposition \ref{spanit}. \eproof
\end{corollary}
\begin{theorem} \label{morextr}
Let $K$ be a cone in $\Rpn$
generated by $S$ and let $u \in S,\
u \neq 0$. Then the following are equivalent:
\begin{enumerate}
\item
$u$ is an extremal in $K$.
\item
For some $j \in \supp (u)$,
$u(j)$ is minimal in $K(j)$.
\item \label{minS}
For some $j \in \supp (u)$, $u(j)$ is minimal in $S(j)$.
\end{enumerate}
\end{theorem}
\begin{proof}
1. $\Longrightarrow$ 3. If $|\supp(u)| = 1$ then $u(j)$ is minimal
in $S(j)$. So suppose that $|\supp(u)| > 1$ and that $u(j)$ is not
minimal in $S(j)$ for any $j \in \supp(u)$. Then for each
$j\in \supp(u)$ there exists $x^j \in S(j)$ such that $x^j \leq
u(j),\ x^j \neq u(j)$. Therefore $u=\bigoplus_{j \in \supp(u)} u_j
x^j$, and $u$ is proportional with none of $x^j$. Hence $u$ is not
an extremal in
$K$.\\
3.  $\Longrightarrow$ 2. Let $v \in K$ and assume that $ j \in
\supp(v)$ and $v(j) \leq u(j)$. We need to show that $v(j) =
u(j)$. By Proposition \ref{spanit}, there is a $w \in S$ such that
$w(j) \leq v(j)$. Thus $w(j) \leq v(j) \leq u(j)$ and by
3. it follows that $w(j) = v(j) = u(j)$.\\
2. $\Longrightarrow$ 1.
Let $u(j)$ be minimal in $K(j)$ for some $j \in \supp(u)$
and suppose that $u = v \oplus w,\ v, w \in K$. Then both $v
\leq u$ and $w \leq u$ and  either  $v_j = u_j $ or $w_j = u_j$,
say (without loss of generality) that  $v_j = u_j$. Hence $v(j)\leq u(j)$
and it follows from 2. that $v(j) = u(j)$. Hence also $v =
u$
which proves 1.\\
\end{proof}
\eproof
Note that in Theorem \ref{morextr} we can of course have
$S = K$. Also note that
Corollary \ref{remove} may be combined with Theorem
\ref{morextr} to yield conditions for a set of
generators to be redundant.
\begin{corollary} \label{genextr}
Let $K$ be a cone in $\Rpn$.
If $D_j(u)$ has a minimal element for each $u \in K$ and each $j
\in \supp(u)$, then $K$ is generated by its extremals.
\end{corollary}
\begin{proof} Suppose that $x^j$ is a minimal element of $D_j(u)$.
Since, for $v \in K(j),\ v \leq x^j$ implies that $v \in D_j(u)$,
$x^j$ is also a minimal element of $K(j)$. We now obtain the
Corollary by combining Proposition \ref{spanit} and Theorem
\ref{morextr}. \end{proof} \eproof
Essentially, the following fundamental result was proved in
\cite[Proposition 2.5.3]{Wag1}. We suitably restate it: every set of
generators $S$ for a cone $K$
can be partitioned as $E\cup F$, where $E$ is a set of
extremals for $K$ and the remainder $F$ is redundant.
Our proof is a combination of
Proposition \ref{spanit} and Theorem \ref{morextr}.
\begin{theorem} \label{mainres}
Let $S$ be a set of scaled generators for a cone $K$ in
$\Rpn$ and let $E$ be the set of scaled extremals in $K$. Then
\begin{enumerate}
\item \label{repeat}$E\subseteq S.$
\item \label{omit}
Let  $F = S\backslash E$. Then for any $u \in F$, the set
$S\backslash\{u\}$ is a set of generators for $K$.
\end{enumerate}
\end{theorem}
\proof
Assertion \ref{repeat} repeats Lemma \ref{extr1}.\\
\msn To prove Assertion \ref{omit}, let $u \in F$. Since $u$ is not an
extremal, by Theorem \ref{morextr} for each $j\in\supp(u)$ there is
$z^j\in K$ such that $z^j(j)<u(j)$. Since $K=\spn(S)$, by
Proposition \ref{spanit} we also have $y^j\in S$ such that
$y^j(j)\le z^j(j)<u(j)$. Evidently $y^j\ne u$, and applying
Proposition \ref{spanit} again, we get that $u$ is a max combination
of $\{y^j : j\in\supp(u)\}$, where $y^j\in S$ are different from $u$.
Thus in any max combination involving $u$, this vector can be replaced
by a max combination of vectors in $S\backslash\{u\}$, and the theorem
is proved.\eproof
The following example shows
that the set $F$ of Theorem \ref{mainres} need not be totally
dependent.
\begin{example} \label{ex1} {\rm Let $K$ be the cone in $R^2_+$
generated by
$u^r = [1,\ 1/r]^T,\ r = 1, \ldots $. The elements
of $K$ scaled with respect to the max norm are $[1,\ a]^T$ with $0<a\le 1$.
Thus $u^1$ is the unique
scaled extremal in $K$. But the set $F = \{u^r:\ r= 2,\ldots\}$ is
not totally dependent since $u^2$ is an extremal in $\spn(F)$
whose scaled elements are $[1,\ a]^T$ with $0<a\le 1/2$.}
\end{example}
The following is a refinement of Theorem \ref{mainres}, and also of
\cite[Theorem 5]{Wag1}.
\begin{theorem} \label{nextres}
Let $E$ be the set of scaled extremals in a max cone $K$. Let $S
\subseteq K$ consist of scaled elements.
Then the following are
equivalent:
\begin{enumerate}
\item
The set $S$ is a minimal set of generators for $K$.
\item
$S = E$ and $S$ generates $K$.
\item
The set $S$ is a basis for $K$.
\end{enumerate}
\end{theorem}
\begin{proof}
1. $\Longrightarrow$ 2. By Theorem \ref{mainres} we have $S = E
\cup F$ where every element of $F$ is redundant in $S$. But since
$S$ is a minimal set of generators, we must have $F = \emptyset$.
Hence $S
=E$.\\
2. $\Longrightarrow$ 3. The set $E$ is independent and a generating set.\\
3. $\Longrightarrow$ 1. By independence of $S$ the span of a
proper subset of $S$ is strictly contained in span($S$).
\end{proof} \eproof
Theorem \ref{nextres} shows that if a cone has a  (scaled) basis
then it must be its set of (scaled) extremals, hence the basis is
essentially unique. We note that a maximal independent set in a cone
$K$ need not be a basis for $K$  as is shown by the following
example.
\begin{example}\label{ex2} {\rm Let $K \subseteq \R^2_+$  consist of
all $[x_1,\ x_2]^T$ with $ x_1 \geq x_2 > 0$. If $1 > a >b > 0$,
then $\{[1,\ a]^T,\ [1,\ b]^T\}$ is a maximal independent set in $K$
which does not generate $K$.}
\end{example}
The following corollary is found e.g. as \cite[Proposition 21]{DS},
\cite[Proposition 2.5]{Jos} and also in \cite[Proposition 1]{Ser},
where it is used to obtain uniqueness results for definite max-plus
matrices. As a special case of this corollary, the standard basis of
$\Rpn$ is essentially unique.
\begin{corollary} \label{finite} If $K$ is a finitely generated cone, then its set
of scaled extremals is the unique scaled basis for $K$.
\end{corollary}
\begin{proof} Since $K$ is finitely generated,
there exists a minimal set of generators
$S$. By Theorem \ref{nextres} $S = E$ and $S$ is a basis.
\end{proof} \eproof
Note that in the tropical geometry \cite{DS} and \cite{Jos} {\em
vertices} of a polytope are defined to be the essentially unique
generators determined in Corollary \ref{finite} and hence vertices
correspond to our extremals (and to Wagneur's irreducible elements).
Next we obtain some corollaries concerning totally dependent sets.
\begin{corollary}\label{totdepinf}
If $S$ is a nonempty scaled totally dependent set in $\Rpn$ then
$S$ is infinite.
\end{corollary}
\begin{proof}
Suppose that $S$ is finite and let $K =$ span($S$). By Corollary
\ref{finite} $K$ contains scaled extremals which, by Theorem
\ref{mainres}, must be contained in $S$ given that $K =$
span($S$). But then $S$ is not totally dependent. This
contradiction proves the result.
\end{proof} \eproof
\begin{corollary}\label{totdep}
Let $K$ be a cone in  $\Rpn$. The following are equivalent:
\begin{enumerate}
\item
There is no extremal in $K$.
\item
There exists a totally dependent set of generators for $K$.
\item
Every set of generators for $K$ is totally dependent.
\end{enumerate}
\end{corollary}
\begin{proof}
Since there always exists a set of generators for $K$ (e.g. $K$
itself), each of the Conditions 2. and 3. is equivalent to Condition
1. by Theorem \ref{mainres}.
\end{proof} \eproof
We now consider $\Rpn$ in the topology induced by the Euclidean
topology of $\Rn$. That is, a set in $\Rpn$ will be called {\em
open} if and only if it is the intersection of an open subset of $\Rn$
with $\Rpn$.
A cone $K$ is called open if
$K\backslash \{0\}$ is open,
%Sergei%
and it is called
closed if it is closed as a subset of $\Rpn$, or
equivalently of $\Rn$.
%Sergei%
\begin{corollary} \label{open} If $K$ is an
open cone in $\Rpn$ that does not contain
unit vectors, then every generating set for $K$ is
totally dependent.
\end{corollary}
\begin{proof} It is enough to show that
there is no extremal in $K$,
for then the result follows
by Theorem \ref{mainres}.
Let $u\in K$. Since $u$ is not a unit vector, there
are at least two indices $k,l\in\supp(u)$.
Since $K$ is open, we have
$w^p=u- \ge e^p \in K,\ p=k,l$ for sufficiently
small $\ge$ and $u=w^k \oplus w^l$.
None of $w^p,\ p=k,l$ is equal to $u$,
hence $u$ is not an extremal,
and the corollary follows.
\end{proof} \eproof
An example of an open cone in $\Rpn$ is furnished by the cone $K$
of all positive vectors in $\Rpn$. We note that, for this
particular case,  Corollary \ref{open} was shown in \cite{CGB}.
Another example of an open cone consists of all vectors
$[a,\ b]^T$ in $\R^2_+$ with $a > b > 0$.
\msn We acknowledge the following Proposition, which is analogous to
Minkowski's Theorem, to \cite[Theorem 3.1]{GK1}. There the result is
proved directly by a minimality argument; here we deduce it from a
corollary to Theorem \ref{morextr} which characterizes extremals of
cones that may not be closed.
It extends
earlier results of \cite{Hel}.
\begin{proposition} \label{c:cc}
Let $K$ be a closed cone in $\Rpn$. Then $K$ is generated
by its set of extremals, and any point in $K$
is a max combination of not more than $n$ extremals.
\end{proposition}
\begin{proof}
Let $u \in K$ and let $j \in \supp(u)$. It is easily shown that
$D_j(u)$ is compact since $K$ is closed. Hence $D_j(u)$ contains a
minimal element $x^j$. The result now follows
by Corollary \ref{genextr} and Corollary \ref{cara}.
\end{proof} \eproof
%Sergei%
The max norm is max linear on $\Rpn$:
%%%%%%%%%
%Old%
%We observe that  max norm $||u||=\max\{u_i\colon i=1,\ldots,n\}$
%on $\Rpn$ is max-linear:%
\begin{equation}
\label{e:linorm} ||\lambda u\oplus\mu v||=\lambda ||u||\oplus\mu
||v||.
\end{equation}
This is exploited in the following proposition.
\begin{proposition}
\label{p:clcl} If $S\subset\Rpn$ is compact and $0\notin S$, then the cone
$K=\spn(S)$ is closed.
\end{proposition}
\begin{proof}
Consider a sequence $u^i\in K$ converging to $v$. Then, by
Corollary \ref{cara} we have
%\begin{equation}
%\label{e:wis}
$$ u^i=\bigoplus_{s=1}^n \lambda_{is} w^{is},$$
%\end{equation}
where $w^{is}\in S$ and $\lambda_{is}\in\R_+$. By~(\ref{e:linorm})
\begin{equation}
\label{e:wisnorm}
||u^i||=\bigoplus_{s=1}^n \lambda_{is}||w^{is}||.
\end{equation}
Since the sequence $u^i$ converges (to $v$), the norms $||u^i||$
are bounded from above by some $M_1>0$. On the other hand,
we have $||w^{is}||\geq M_2$ for some $M_2>0$, since $S$ is closed and does
not contain $0$. Then by~(\ref{e:wisnorm})
$\lambda_{is}||w^{is}||\leq M_1$ for all $i$ and $s$,
and $\lambda_{is}\leq M_1 M_2^{-1}$ for
all $i$ and $s$. Thus $\lambda_{is}$ are bounded from above. But
$||w^{is}||$ are also bounded from above, since $S$ is compact.
This implies that there is a subsequence $u^{j(i)}$ such that for
all $s=1,\ldots,n$ the sequences $w^{j(i)s}$ and $\lambda_{j(i)s}$
converge. Denote their limits by $\Bar{w}^s$ and
$\Bar{\lambda}_s$, respectively, then $\Bar{w}^s\in S$ and
$\Bar{\lambda}_s\leq M_1 M_2^{-1}$. By continuity of $\oplus$ and
$\otimes$ we obtain that
\begin{equation}
\label{e:barws}
v=\bigoplus_{s=1}^n \Bar{\lambda}_s \Bar{w}^s.
\end{equation}
\end{proof}
Thus $v\in K$.\eproof
\begin{corollary}
\label{c:extr} If the set of scaled extremals of a max cone $K$ is
closed and generates $K$, then $K$ is closed.
\end{corollary}
\begin{corollary}
\label{c:fincl}
Any finitely generated max cone $K$ is closed.
\end{corollary}
We now give a counterexample to the converses of Corollary~\ref{genextr}
and Proposition~\ref{c:cc} (part 1), and to the converse of
Corollary~\ref{c:extr} (part 2).
\begin{example}\label{counterex} {\rm
\begin{enumerate}
\item
In $\R^3_+$ let $S$ consist of all vectors $[x_1,\ x_2,\ 1]^T,\
0\leq x_1 < 1/2$ such that $x_1 + x_2 =1 $ and let $K = \spn(S)$.
Then the section of $K$ given by $x_3=1$ 
consists of all vectors 
$[x_1,\ x_2,\ 1]^T,\ 0\leq x_1 < 1/2,\ 0 \leq x_2 \leq 1 $ such that
$x_1 + x_2 \geq 1 $. Note that $S$ is the set of extremals of $K$
scaled with respect to the max norm, but $K=\spn(S)$ is not closed
and for any $u\in K$ there are no minimal elements in $D_1(u)$ and
$D_2(u)$.
\item
Now let $S' = S \cup \{u\}$, where $u = [1/2,\ 0,\ 1]^T$ and let $K'
= \spn(S')$. Then the section of $K'$ given by $x_3 = 1$  consists
of $K$ together with the line segment whose end points are $u$ and
$[1/2,\ 1,\ 1]^T$. Thus $K'$ is closed. The set of scaled extremals
of $K'$ is $S'$ which is not closed.
\end{enumerate} }
\end{example}
The cross sections of $K$ and $K'$ by $x_3=1$ are shown on Figure 1, together
with the generating sets $S$ and $S'=S\cup\{u\}$.

%\begin{figure}[h]
%\centering
%\vskip-1cm
%\includegraphics[width=10cm]{counterex}
%\vskip-7cm
%\caption{Max cones of Example \ref{counterex}}
%\end{figure}

\begin{figure}[htbp]
%\vskip0.5cm
\hskip1cm
\epsfysize=5cm
\epsffile{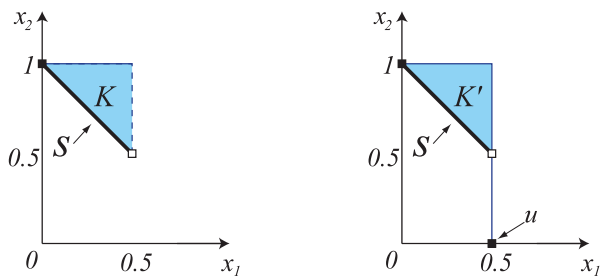}
\vskip -0.5cm
\caption{Max cones of Example \ref{counterex}}
\end{figure}

\section{Algorithmic considerations}
As explained in the introduction we redefine our basic concepts for
this section which is concerned with finitely generated cones. We
also restate a suitable adaptation of Corollary \ref{finite}.
\begin{definition}{\rm
Let $V \in \Rpnk$ and let $V_{\widehat{i}}$ be the matrix obtained
from $V$ by deleting column $i,\ i = 1,\ldots,k$. Then the cone
$K$ {\em generated} by the columns $v^1,\ldots,v^k$ of $V$
consists of all vectors of form $V\otimes x,\ x \in \R_+^k$.
Further, the columns of $V$ form a {\em basis} for $K$ if, for
$i,\ i = 1,\ldots,k$, there is no $x \in \R_+^{k-1}$ such that
$V_{\widehat{i}} \otimes x = v^i$.}
\end{definition}
\begin{proposition}
Let $V \in \Rpnk$. Then there exists a submatrix $U \in
\R^{np}_+, 0 \leq p \leq k$ whose columns form a basis for the
cone generated by the columns of $V$ (and every other basis is of
form $UPD$, where $P$ is a permutation matrix and  $D$ is a
diagonal matrix with nonzero diagonal elements). \eproof
\end{proposition}
We shall apply the following proposition. Note that all statements
in this proposition have been proved in a more general setting in
\cite{CG1}. See also \cite{But}, \cite[Chapter III]{CG3} and
\cite{Vor}.
\begin{proposition} \label{extension} Let $U \in \Rpnk$ with all columns nonzero
and let $v \in \Rpn$. Let $x \in \Rpk$ be defined by
\begin{equation} \label{critical}
x_i = \min\{v_j/u^i_j :\ u^i_j \neq 0, j = 1, \ldots, n\}
\end{equation}
for $i = 1,\ldots,k$. Then
\begin{eqnarray}
\label{e1} U \otimes x & \leq & v, \\
\label{e2}  x & =& \max\{z \in \Rpk : U \otimes z \leq v\},\\
\label{e3} U \otimes x & = & \max\{U \otimes z: z \in \Rpk, U
\otimes z \leq v\}.
%\label{e4} \exists z \in \Rpk\quad, U \otimes
%z = v & \Longleftrightarrow & U \otimes x = v.
\end{eqnarray}
Further, there exists $z \in \Rpk$ such that $U \otimes z = v$ if
and only if $U \otimes x = v$.
%\end{enumerate}
\end{proposition}
\proof Assertion (\ref{e1}) follows from the observation that $U
\otimes z \leq v$ if and only if $z_i \leq v_j/u^i_j$ if $j \in
\supp(u^i),\ i = 1,\ldots,k$. Note that $x \in \Rpk$ since no
column of $U$ is zero. Since $\otimes $ is isotone (that is,
$x\leq y$ implies $A\otimes x\leq A\otimes y$), assertions
(\ref{e2}) and (\ref{e3}) follow immediately. For the final
statement assume that $U\otimes z=v$ for some $z$. By (\ref{e1})
and (\ref{e3}) we have $ v=U\otimes z\leq U\otimes x\leq v$, and
the statement follows. The converse is trivial. \eproof
%Sergei%
\begin{algorithm} \label{algo} {\rm
{\em Input:} $V \in \Rpnk$.\\
{\em Output:} An $n \times p$ submatrix $U$ of $V$ whose columns
form the essentially unique basis for the cone generated by the
columns of $V$.
\medskip\noindent\\
{\em Step 1}. Initialize $U = V$.\\
{\em Step 2}. For each $j = 1,\ldots,k$ if
$u^j\neq 0$ set
$v = u^j$, and for each $i\neq j$ compute $x_i$ by
(\ref{critical}), if $u^i\neq 0$, and set $x_i=0$
otherwise. If $U_{\widehat{j}}x=v$, set $u^j=0$.\\
{\em Step 3}. Delete the zero columns of $U$. The remaining
columns of $U$ are the basis we seek.}
\end{algorithm}
\begin{remark}{\rm
The restriction in Proposition \ref{extension} that each column $U
\in \Rpnk$ must have a positive element was imposed to avoid
definitions for $a/0, a > 0$, or $0/0$. The restriction is
inessential in the sense that for general $U \in \Rpnk$ we may
define $x_i$ by (\ref{critical}) whenever $u^i \neq 0$ and choose
$x_i$ arbitrarily in $\Rpk$ whenever $u^i = 0$. Then all
assertions of the Proposition still hold, with exception of
(\ref{e2}). It is possible to extend $\Rpn$ by adding a maximal
element $\infty$ so that (\ref{e2}) still holds.}
\end{remark}
We omit details and present the MATLAB program maxbas that
implements Algorithm \ref{algo} but employs such an extension. We
also give an example with some
elements equal to $0$. Note that in \cite[Theorem 16.2]{CG2} a
related algorithm called {\em A-test} has been presented. It
enables us to identify columns that are dependent on other columns
of an $n\times k$ matrix in $O(nk^2)$ time. However, there is no
discussion of bases in connection with this method in \cite{CG2}.
\begin{verbatim}
%the unique max times basis for the max col space of A
%function [B,f] = maxbas(A),
%B = the unique max times basis for the max col space of A
%f = indices of columns of B in A
%calls maxpr,  max multiplication  of matrices
 function  [B,f] = maxbas(A)
 [m,n] = size(A); B = A; t = max(max(A));
 for j  = 1:n
     v = compl(j,n);
     c = B(:,j); BB = B(:,v); warning('off'),
     e = ones(1,n-1); C = c*e; x = min(C./BB)';
     z= maxpr(BB,x);
     if abs(c-z) < t*eps, B(:,j) = 0; end,
end u = max(B); f = find(u >t*eps); B = B(:,f);
A =
  Columns 1 through 5
     1     9    10     5     9
     2    10    10     0    10
     3    15    14     7     0
     4    20    16     8    12
>> [B,f] = maxbas(A)
B =
     1     5     9
     2     0    10
     3     7     0
     4     8    12
f =
     1     4     5
\end{verbatim}
We note that a second form of the algorithm may be based on
\emph{set covering} condition~(\ref{setcov}) below,
which appears in~\cite[Theorem 15.6]{CG2}. It can also
be found in \cite{Vor} Section 2 and in
\cite{But} Section 2 but
only in the case when all vectors are positive. See also
\cite[Theorem 3.5]{AGK}
for an interesting functional generalization of this condition (and more).
%It may also be easily derived from Proposition
%\ref{spanit}.
With $v$ and $U$ as in Step 2, denote by $N_i$ the
set $\{j\colon v(j)\geq u^i(j)\}$.
By Proposition \ref{spanit},
$v\in\spn(u^1,\ldots,u^m)$ if and only if
\begin{equation}
\label{setcov} \bigcup_{i=1}^m N_i=\supp(v).
\end{equation}
With $x$ given by
(\ref{critical}), we note that
\begin{equation}
\label{nidef}
N_i=
\begin{cases}
%\text{arg}\min_\{v_j(u^i_j)^{-1}\colon j\in\supp(u^i)\},&\text{if $\Bar{x}_i\ne 0$};\\
\{j\in \supp(u^i):  v_j/u^i_j = x_i\}  &  {\rm if}\ x_i \neq 0,\\
\emptyset,&\text{if $x_i=0$.}
\end{cases}
\end{equation}
Thus Step 2 in Algorithm \ref{algo} may be replaced by

\msn {\em Step 2'}:  For each $j=1,\ldots,k$ such that $u^j\ne 0$:
set $v=u^j$ and for each $i\ne j$ compute $N_i=\{j\colon v(j)\geq
u^i(j)\}$ according to~(\ref{nidef}), if $u^i\ne 0$, and set
$N_i=\emptyset$ otherwise. If $\bigcup_{i\neq j} N_i=\supp(v)$,
set $u^j=0$.

\msn The version with Step 2' is also well-known. It is implemented
in the max-plus toolbox of Scilab, a freely distributed software. See
\cite[Sect.III-B]{GS} for the documentation.

\msn Our algorithms are of complexity $O(nk^2)$.

\msn If $S$ is the set of columns of the matrix $U$, then it
follows from Theorem \ref{morextr} that a basis for the cone
generated by $S$ consists of the union of the $n$ sets $M(j), \ j
= 1,\ldots,n$, where $M(j)$ consist of the vectors minimal in
$S(j)$.  The problem of finding all maxima (or minima) of $k$ vectors in
$\Rn$ is considered in \cite{KLP}, and also in
\cite[Section 4.1.3]{PS}, where it is dubbed the
problem of Erehwon Kings. The computational
complexity of methods developed in \cite{KLP}
and \cite{PS} is bounded from above by $O(n^2 k(\log_2 k)^{n-2})+
O(k\log_2 k),\,n\geq 2$, see
\cite[Theorem 4.9]{PS} and \cite[Theorem 5.2]{KLP}.
\footnote{That is, the complexity is not greater than
$O(k\log_2 k)$ if $n=2$
and not greater than $O(n^2 k(\log_2 k)^{n-2})$ if $n\geq 3$.}
To solve our problem
we can apply these methods
to each $S(j),j=1,\ldots,n$ separately.
Taking into account that for
each $j$ we need $O(nk)$ operations to find the coordinates of essentially
$(n-1)$-dimensional vectors in $S(j)$, this yields an alternative method with
complexity not smaller than $O(n^2 k)$ and not greater than
$O(n^3 k(\log_2 k)^{n-3})+O(k\log_2 k),\, n\geq 3$.
This method may be preferred if $\log_2 k$ is substantially larger than $n$.
%Its instance for $n=3$ is the Triple Sorting Algorithm described in
%\cite[Section 5]{Jos}.

\msn \textbf{Acknowledgment.} We thank G.M. Engel and an unknown
referee for valuable comments. We thank S.~Gaubert and R.~Katz for
informing us of their papers \cite{GK0} and \cite{GK1} and for very
useful discussion. All have helped to improve this paper.

\end{document}